\documentclass[11pt,reqno]{amsart}
\usepackage{amsfonts,amssymb}
\usepackage{xy,epsfig}

\catcode `\@=11

\def\numberbysection{\@addtoreset{equation}{section}
\renewcommand{\theequation}{\thesection.\arabic{equation}}}
\numberbysection
\def\subsubsection{\@startsection{subsubsection}{3}%
  \normalparindent{.5\linespacing\@plus.7\linespacing}{-.5em}%
  {\normalfont\bfseries}}
\def\mapright#1{\smash{\mathop{\lra}\limits^{#1}}}

\def\lra{\longrightarrow}
\def\rto{\rightarrow}

\def\S{{\mathbb S}}

\def\PPP{{\mathbb P}}
\def\PP{{\mathcal P}}
\def\I{{\mathbb I}}
\def\II{{\mathcal I}}

\DeclareMathOperator{\st}{st}
\DeclareMathOperator{\inje}{inj}
\DeclareMathOperator{\Id}{Id}
\DeclareMathOperator{\Ker}{Ker}

\DeclareMathOperator{\Prim}{Prim}
\DeclareMathOperator{\Sh}{Sh}

\DeclareMathOperator{\Com}{\mathcal{C}om}
\DeclareMathOperator{\As}{\mathcal{A}s}
\DeclareMathOperator{\Lie}{\mathcal{L}ie}
\DeclareMathOperator{\Fin}{\mathbf{Fin}}
\DeclareMathOperator{\Inj}{\mathbf{Inj}}
\DeclareMathOperator{\Vect}{\mathbf{Vect}}
\DeclareMathOperator{\Coalg}{\mathbf{Coalg}}
\DeclareMathOperator{\kk}{\mathbf{k}}
\DeclareMathOperator{\un}{\mathbf{1}}
\DeclareMathOperator{\Pois}{\mathcal{P}ois}
\DeclareMathOperator{\Mag}{\mathcal{M}ag}

\xyoption{all}
\CompileMatrices

\begin{document}
\title[Lie theory for Hopf operads]{Lie theory for Hopf operads}
\author[M. Livernet]{Muriel Livernet}
\address{ML: Institut Galil\'ee\\
Universit\'e Paris Nord\\
  93430 Villetaneuse       \\
         France}
\email{livernet@math.univ-paris13.fr}
\urladdr{http://www.math.univ-paris13.fr/$\sim${}livernet/}

\author[F. Patras]{Fr\'ed\'eric Patras}
\address{FP: CNRS et Universit\'e de Nice-Sophia Antipolis \\
Laboratoire J.-A. Dieudonn\'e \\
Parc Valrose\\
06108 Nice cedex 02  France}
\email{patras@math.unice.fr}
\urladdr{http://math.unice.fr/$\sim${}patras/}

\thanks{The first-named author thanks the Institut Mittag-Leffler (Djursholm, Sweden) 
where a part of this work was carried out during her stay. The second-named author 
was supported by the ANR grant AHBE 05-42234}

\date{\today}

\begin{abstract}
The present article takes advantage of the properties of algebras in the category of 
$\S$-modules (twisted algebras) to investigate further the fine algebraic structure of Hopf operads. 
We prove that any Hopf operad $\PP$ carries naturally the structure of twisted Hopf $\PP$-algebra. 
Many properties of classical Hopf algebraic structures are then shown to be encapsulated in the twisted 
Hopf algebraic structure of the corresponding Hopf operad. In particular, various classical 
theorems of Lie theory relating Lie polynomials to words (i.e. elements of the tensor algebra) are 
lifted to arbitrary Hopf operads. 
\end{abstract}

\maketitle

\section*{Introduction}
Let $\PP$ be an arbitrary algebraic operad, that is, a monoid in the category of 
$\S$-modules or, equivalently, the analytic functor associated to a given (suitable) class of 
algebras. One can extend the usual definition of algebras over $\PP$ in the category of vector 
spaces and define algebras over $\PP$ in
the category of $\S$-modules. These algebras are classically refered to as \it 
twisted $\PP$-algebras\rm . 
The meaning and usefulness of twisted algebraic structures was pointed
out already in 1978 by Barratt, whose work was motivated by the study
of homotopy invariants \cite{Barratt78}. He introduced the notion of
\it twisted \rm Lie algebras (Lie algebras in the category of
$\S$-modules).
His constructions were later extended to various other algebraic structures by Joyal \cite{Joyal86}.

The general definition of an algebra over an operad $\PP$ in the
category of $\S$-modules is more recent. It appears for example in the
work
of Fresse \cite{Fr04} under the categorical name of ``left
$\PP$-modules''. We refer to his article, also for further
bibliographical informations on operads.

The purpose of the present article is to take advantage of the
properties of twisted algebraic structures to investigate further
the internal algebraic structures of operads. We are mainly interested
in Hopf operads, that is, the particular class of operads $\PP$ for
which one can provide the tensor product of two $\PP$-algebras with the structure of a $\PP$-algebra.

We first show that a Hopf operad $\PP$ has naturally a twisted Hopf
$\PP$-algebra structure. We then show that set of primitive elements for this
twisted Hopf $\PP$-algebra structure is a sub-operad of $\PP$. For
example, when $\PP =\As$, the associative operad, there is a twisted
Hopf algebra structure on the direct sum $S_\ast$ of the symmetric
group algebras.   As could be expected, its operad of
primitive elements is the Lie operad, which makes more precise the results obtained in \cite{PatReu04} on the fine twisted Hopf algebra structure of $S_\ast$.

We then turn to reciprocity laws, namely to the relations between the
internal algebraic structure of an operad and the structure of the
algebras
over this operad. More precisely, we are interested in the link
between the twisted Hopf $\PP$-algebra structure of the Hopf operad
$\PP$,   its primitive suboperad and
the structural properties of  Hopf $\PP$-algebras. We show that many
classical properties in the theory of free associative algebras and
free Lie algebras go over to arbitrary Hopf operads and their primitive suboperads.

Our general results are illustrated on the most classical examples of Hopf operads, namely the associative and Poisson operads.

%%%%%% SECTION: TWISTED ALGEBRAS OVER AN OPERAD

\section{twisted algebras over an operad}

For this section we refer to May \cite{May72}, Barratt \cite{Barratt78}, Joyal \cite{Joyal86},
Fresse \cite{Fr04} and Patras-Reutenauer \cite{PatReu04}.

%%%%%%%%%% SOUS-SECTION:  Les S-modules

\subsection{$\S$-modules and some related categories of modules.}

%%%%%%%%%%%%% les S-modules.

\subsubsection{$\S$-modules.}
A {\sl $\S$-module} $M=\{M(n)\}_{n\geq 0}$
is a collection of right $S_n$-modules over a ground field $\kk$.
Notice that we will write permutations in $S_n$
as sequences: $\sigma =(\sigma (1),...,\sigma (n))$.

Following Joyal, a $\S$-module is
equivalent to a {\sl vector species} that is a contravariant functor from the category of
finite sets $\Fin$ (and set isomorphisms) to the category $\Vect$ of vector
spaces over a field $\kk$. This equivalence goes as follows:
from a $\S$-module $M$ one defines the vector species
$$S\mapsto {{M}}(S):=\bigoplus_{i_*:S\rto \{1,\ldots,r\}} M(r)/\equiv$$
where $r$ is the number of elements of $S$ and $i_*$ is a bijection.
The equivalence relation is given by $(m\cdot\sigma,i_*)\equiv
(m,\sigma\circ i_*)$.

Conversely, the skeleton of a vector
species ${M}$ is a $\S$-module. The action of $S_n$ on $M(n):={{M}}(\{1,\ldots,n\})$ is given by
 $m\cdot\sigma ={M}(\sigma)(m)$.

%%%%%%%%%%%%%%%%%%%%%%%Les I-modules

\subsubsection{$\I$-modules}\label{S-Imodules}
As it is well-known, and recalled in the next section,
$\S$-modules or, equivalently, vector species, allow to study algebraic structures
from the systematical point of view given by the operadic framework. As we will
be interested with algebraic structures provided with additional properties, such as
the existence of a unit, let us introduce a generalization of $\S$-modules suited for
our purposes.

Let $\Inj$ be the category of finite sets and injections. An \sl injective species \rm is a
contravariant functor from $\Inj$ to vector spaces. Correspondingly, a
$\I$-module $J=\{J(n)\}_{n\geq 0}$ is a collection of right $S_n$-modules together
with \sl degeneracy maps \rm 

$$\partial_i :J(n)\longmapsto J(n-1), i=1,...,n$$
such that: 
\begin{itemize}
\item For any $1\leq i\leq n$, $1\leq j\leq n-1$, with $j\geq i$, we have the equality 
between maps from $J(n)$ to $J(n-2)$:
\begin{equation}\label{E-degeneracy1}
\partial_j\circ\partial_i =\partial_i\circ\partial_{j+1}.
\end{equation}
\item For any $1\leq i\leq n$ and any $m\in J(n),\sigma\in S_n$, we have:
\begin{equation}\label{E-degeneracy2}
\partial_i(m\cdot\sigma) =\partial_{\sigma (i)}(m)\cdot \partial_i(\sigma),
\end{equation}
with $\partial_i(\sigma)=\st(\sigma(1),...\sigma(i-1),\sigma(i+1),...,\sigma(n))$ where $\st$ 
states for the
standardization of a sequence of distinct integers.
\end{itemize}

Recall that, in general, the standardization of a sequence of lenght $p$  
of distinct nonnegative integers is the process by which the elements of the sequence are 
replaced by the integers $1,...,p$ in such a way that the relative order of the elements 
in the sequence is preserved. The process is better understood by means of an example: 
if $\sigma =(3,2,6,1,8,7,5,4)\in S_8$, we have
$$\st(\sigma (1),\sigma(2),\sigma(4),...,\sigma (8))=\st(3,2,1,8,7,5,4)=(3,2,1,7,6,5,4)\in S_7.$$

The relation between the two notions of $\I$-modules and injective species is given by 
the same formula as the equivalence between the two notions of $\S$-modules and vector 
species. The map $\partial_i$ from $J(n)$ to $J(n-1)$ corresponds to the \it elementary 
injection \rm $\inje_i$ from $[n-1]$ to $[n]$ defined by: 
$$\inje_i(k):=\begin{cases} k &  {\rm if}\ k\leq i-1 \\ 
k+1 &  {\rm if}\ k\geq i.\end{cases}$$ 

The equivalence follows from the observation that any injection 
from $[k]$ to $[n]$, $k<n$, factorizes uniquely as a composition of permutations and  elementary injections:
$$\inje_{i_{n-k}}\circ ...\circ \inje_{i_1}\circ \sigma$$ with $\sigma\in S_k$ and $i_1<...<i_{n-k}$.

Notice that similar notions have been studied by Cohen, May and Taylor in \cite{CMT78}   
and Berger \cite{Berger96}.

%%%%%%%%%%%%%%%%%%%%%exemple du groupe symetrique

\subsubsection{The symmetric groups as an $\I$-module.}\label{S-symmetricgroups}
Consider now the $\S$-module given by $\kk[S_n]$ for all $n\geq 0$ and the
right action given by right multiplication.
Then $\partial_i(\sigma)=\st(\sigma(1),...\sigma(i-1),\sigma(i+1),...,\sigma(n))$ satisfies
relations (\ref{E-degeneracy1}) and (\ref{E-degeneracy2}), thus this $\S$-module can be 
provided with the structure of an $\I$-module.

  The property also follows from the observation that $\S$ can be viewed, from the categorical 
point of view, as the skeleton of the full subcategory $\Fin$ of $\Inj$, whereas $\I$ identifies, 
categorically, with the skeleton of $\Inj$. In this point of view, $\S$-modules and $\I$-modules 
identify with contravariant functors from $\S$ and $\I$.
%As a consequence,
%if equation (\ref{E-degeneracy1}) is satisfied for a
%system of generators of the $\S$-module $J$ then imposing equation (\ref{E-degeneracy2})
%implies that equation (\ref{E-degeneracy1}) is fullfilled in $J$.
%This observation will be useful in order to check that a $\S$-module is also a $\I$-module.

%%%%%%%%%%%%%%%%%%% Operations

\subsubsection{Operations on $\S$-modules and $\I$-modules.} \hfill\break

\noindent{i)} The categories of vector species and injective species are linear symmetric monoidal categories with the
following {\sl tensor product}:
$$({{M}}\hat\otimes {{N}})(U)=\bigoplus_{I\sqcup J=U} {M}(I)\otimes {N}(J)$$
where $I\sqcup J$ runs over the partitions of $U$.

Explicitely, a map $\phi$ from $V$ to $U$ in $\Inj$ induces a map from
$(M\hat\otimes N)(U)$ to $(M\hat\otimes N)(V)$. 
Its restriction to ${M}(I)\otimes {N}(J)$, where $I\sqcup J=U$, is the
tensor product of the maps from $M(I)$ to $M(I')$ and from 
$N(J)$ to $N(J')$ that are induced by the restriction of $\phi$ to a map from $I':=\phi^{-1}(I)$ to $I$ (resp. from $J':=\phi^{-1}(J)$ to $J$).

Translated to
$\S$-modules and $\I$-modules, it gives
$$(M\hat\otimes N)(n)=\bigoplus_{I\sqcup J=[n]} M(I)\otimes M(J).$$ This definition makes use of
the equivalence between $\S$-modules and vector species and between $\I$-modules and
injective species. The equivalent
definition using $\S$-modules and $\I$-modules only is the following
\begin{align*}
(M\hat\otimes N)(n)=&\bigoplus_{p+q=n} (M(p)\otimes N(q))\otimes_{S_p\times S_q} \kk[S_n]\\
=&\bigoplus_{p+q=n} (M(p)\otimes M(q))\otimes\kk[\Sh_{p,q}]
\end{align*}
where $\Sh_ {p,q}$ is the set of $(p,q)$-shuffles that is permutations of $S_n$
written $(\tau_1,\ldots,\tau_p,\rho_1,\ldots,\rho_q)^{-1}$
with $\tau_1<\ldots<\tau_p$ and $\rho_1<\ldots<\rho_q$. The second
equality is a consequence of the unique decomposition of any permutation $\sigma\in S_n$ as
$\sigma=(\sigma_1\times\sigma_2)\cdot \alpha$ where $\alpha$ is a
$(p,q)$-shuffle. The unit for this tensor product is the 
$\S$-module (resp. $\I$-module)
$\un$ given by
$$\un(n)=\begin{cases} \kk, & {\rm if}\ n=0,\\ 0,& \ {\rm otherwise}.
\end{cases}$$

\medskip

When $M$ and $N$ are $\I$-modules, the structure of $\I$-module on $M\hat\otimes N$ is
given as follows: for any $m\in M(p)$, any $n\in N(q)$ with $p+q=n$,
\begin{equation}\label{E-partialontensor}
\partial_i(m\otimes n)=\begin{cases} \partial_i(m)\otimes n & {\rm if}\ 1\leq i\leq p, \\
m\otimes \partial_{i-p}(n) & {\rm if}\ p+1\leq i\leq p+q.\end{cases}
\end{equation}
Notice that this relation determines entirely the  $\I$-module structure on 
$M\hat\otimes N$ thanks to
equation (\ref{E-degeneracy2}), since $M(p)\otimes N(q)$ generates 
$(M\hat\otimes N)(p+q)$ as an $S_{p+q}$-module.

\medskip

The symmetry isomorphism $\tau_{M,N}:M\hat\otimes N\rto N\hat\otimes M$ is given by
$$\tau_{{M},{N}}: {{M}}(I)\otimes {{N}}(J) \rto {N}(J)\otimes {M}(I)$$
in the species framework,
that is, in the $\S$ and $\I$-modules framework
\begin{equation}\label{E-symmetryiso}
\tau_{M,N}(m\otimes n)=(n\otimes m)\cdot \zeta_{p,q}
\end{equation} 
where $m\in M(p), n\in N(q)$ and
$\zeta_{p,q}=(q+1,\ldots,q+p,1,\ldots,q)$. 
The symmetry isomorphism is a morphism of $\I$-modules that is
\begin{equation}\label{E-Isymmetryiso}
\partial_i\tau_{M,N}=\tau_{M,N}\partial_i.
\end{equation} 
Indeed, if $1\leq i\leq p$ one has
\begin{multline*}
\partial_i\tau_{M,N}(m\otimes n)=\partial_i((n\otimes m)\cdot \zeta_{p,q})
=\partial_{\zeta_{p,q}(i)}(n\otimes m)\cdot\partial_i(\zeta_{p,q}) \\
=\partial_{q+i}(n\otimes m)\cdot \zeta_{p-1,q}=(n\otimes \partial_i(m))\cdot \zeta_{p-1,q}
=\tau_{M,N}(\partial_i(m)\otimes n)\\
=\tau_{M,N}\partial_i(m\otimes n).
\end{multline*}
The same computation holds if $p+1\leq i\leq p+q$.

\medskip

For any $\sigma\in S_n$, the symmetry isomorphism induces an isomorphism 
$\tau_\sigma$ of $\S$-modules and $\I$-modules from
 $M_1\hat\otimes\ldots\hat\otimes M_n$ to  $M_{\sigma^{-1}(1)}
\hat\otimes\ldots\hat\otimes M_{\sigma^{-1}(n)}$ given by
\begin{equation}\label{F-symmetry}
\tau_\sigma(m_1\otimes\ldots\otimes
m_n)=(m_{\sigma^{-1}(1)}\otimes\ldots\otimes m_{\sigma^{-1}(n)}) \cdot
\sigma(l_1,\ldots,l_n)
\end{equation}
where $m_i\in M_i(l_i)$ and
 $\sigma(l_1,\ldots,l_k)$ is the permutation of $S_{l_1+...+l_k}$ obtained by 
replacing $\sigma(i)$ by the block $\Id_{l_i}$. More precisely
$$\sigma(l_1,\ldots,l_k)=(B_1,\ldots,B_k)$$
where $B_i$ is the block
$l_{\sigma^{-1}(1)}+\ldots+l_{\sigma^{-1}(\sigma(i)-1)}+[l_i].$
For instance 
$$(2,3,1)(a,b,c)=(c+1,\ldots,c+a,c+a+1,\ldots,c+a+b,1,\ldots,c).$$

\medskip

\noindent ii) The notation $M\otimes N$ is devoted to the $\S$-module or 
$\I$-module given by the collection
$M(n)\otimes N(n)$ together with the diagonal action of maps. It is the vector or 
injective species $({M}\otimes {N})(U)={M}(U)\otimes {N}(U)$.

\medskip

\noindent iii)  The category of $\S$-modules is endowed with another 
monoidal structure (which is not symmetric): the {\sl plethysm $\circ$} defined by
$$(M\circ N)(n):=\bigoplus_{k\geq 0} M(k)\otimes_{S_k} (N^{\hat\otimes k})(n),$$
where $S_k$ acts on the left of $(N^{\hat\otimes k})$ by formula (\ref{F-symmetry}).

An alternative definition using vector species is given by
$$({M}\circ {N})(U):= \bigoplus_{k\geq 1} 
M(k)\otimes_{S_k}(\bigoplus_{I_1\sqcup\ldots\sqcup I_k=U} 
{N}(I_1)\otimes\ldots\otimes {N}(I_k)).$$
when $U\not=\emptyset$ and
where the action of $S_k$ is given by
$$\sigma\cdot ( {N}(I_1)\otimes\ldots\otimes {N}(I_k))={N}(I_{\sigma^{
-1}(1)})\otimes\ldots\otimes {N}(I_{\sigma^{-1}(k)}).$$

The unit for the plethysm is the $\S$-module $I$ given by
$$I(n)=\begin{cases} \kk, & {\rm if}\ n=1,\\ 0,& \ {\rm otherwise}.
\end{cases}$$

\medskip

Because of relation (\ref{E-Isymmetryiso}), the left action of $S_k$ on
$N^{\hat\otimes k}$ commutes with the action of the morphisms in $\Inj$.  
Therefore, one can also define the plethystic product $M\circ {N}$ of a $\S$-module $M$ with a
$\I$-module $N$ by the same formula
$$(M\circ {N})(U):= \bigoplus_{k\geq 0}
M(k)\otimes_{S_k} N^{\hat\otimes k}(U),$$
together with a $\I$-module structure defined
for any $m\in M(k), \bar n\in N^{\hat\otimes k}(l)$, by:
$$\partial_i(m\otimes \bar n):=m\otimes \partial_i(\bar n).$$

In particular, since any $\I$-module is also a $\S$-module, the plethystic product of
two injective species is an injective species. However some properties that hold
for vector species break down for injective species. For example, the $\S$-module
$I$ can be provided  with a unique $\I$-module structure and is still a left unit for
the plethystic product of $\I$-modules. However, it is not a right unit: the
$\I$-module structure of $M\circ I$ is not the $\I$-module structure of $M$, except in
the very particular case when the degeneracy maps from $M(n)$ to $M(n-1)$ are all zero.

%%%%%%%%%%%%%%% transformations naturelles

\subsubsection{Natural transformations relating the operations on $\S$ and $\I$-modules}

One has two natural transformations in $A,B,C$ and $D$:

$$\begin{array}{cccc}
T_1:& (A\otimes B)\circ (C\otimes D)&\rto& (A\circ C)\otimes (B\circ D)\\
T_2:& (A\otimes B)\circ (C\hat\otimes D)&\rto& (A\circ C)\hat\otimes (B\circ D)
\end{array}$$
obtained by interverting terms in the   corresponding direct sums. The symmetry
isomorphism $\tau$ has to be taken into account in order to define $T_2$.
Since $\tau$ commutes  with degeneracies,  $T_2$ is a morphism of $\I$-modules if $C$ and $D$ are $\I$-modules.

In terms of $\S$-modules, one has the following description: let $a\in A(k)$ and
$b\in B(k)$, $c_i\in C(l_i),
d_i\in D(l_i)$, $e_i\in C(r_i)$, $f_i\in D(s_i)$:

\begin{multline*}
T_1((a\otimes b) \otimes (c_1,d_1,c_2,d_2,\ldots,c_k,d_k))=(a\otimes
c_1,\ldots,c_k) \otimes (b\otimes d_1,\ldots d_k) \hfill \\
\ \ \ T_2((a\otimes b) \otimes (e_1,f_1,e_2,f_2,\ldots,e_k,f_k))= \hfill \\
(a\otimes e_1,\ldots,e_k) \otimes (b\otimes f_1,\ldots f_k)\cdot
\sigma(r_1,s_1,r_2,s_2,\ldots,r_k,s_k)
\end{multline*}
where $\sigma=(1,k+1,2,k+2,3,\ldots,2k-1,2k)$.

The relations are even more natural when written in terms of vector species.
Let us consider, for example, the relation defining $T_2$. We have, for
$e_i\in C(S_i)$ and $f_i\in D(T_i)$ with $\coprod_{i=1}^k(S_i \coprod T_i)=U$:

\begin{multline*}
T_2: (A(k)\otimes B(k))\otimes_{S_k}((C(S_1)\otimes D(T_1))\otimes ...\otimes
(C(S_k)\otimes D(T_k)))\\
\longrightarrow (A(k)\otimes_{S_k}( C(S_1)\otimes ...\otimes C(S_k)))\otimes (B(k)\otimes_{S_k}(D(T_1)\otimes ...\otimes D(T_k)))
\end{multline*}
$$T_2(a\otimes b \otimes (e_1,f_1,e_2,f_2,\ldots,e_k,f_k))=(a\otimes e_1,\ldots,e_k) \otimes (b\otimes f_1,\ldots f_k).$$

  As a consequence of the definitions, the following relations hold
\begin{align}\label{R-iso}
(A\hat\otimes B)\circ C=&(A\circ C)\hat\otimes (B\circ C),\\
\label{R_T1_T1}
T_1(T_1\circ \Id)=&T_1(\Id\circ T_1),\\
\label{R_T1_T2}
T_2(T_1\circ \Id)=&T_2(\Id\circ T_2).
\end{align}

Notice that the last identity expresses in two different ways the
natural transformation
$$(A\otimes B)\circ (C\otimes D)\circ (E\hat\otimes F)\rto (A\circ C\circ E)\hat\otimes (B\circ D\circ F).$$

%%%%%%%%%%%%%%%%%%%%SOUS-SECTION : LES OPERADES

\subsection{Operads}

%%%%%%%%%%%%%%% def: operad

\subsubsection{Definition}\label{D-operad} An {\sl operad} is a monoid in the category of $\S$-modules
with respect to the plethysm. Hence an operad is a $\S$-module $\PP$ together with a product $\mu_\PP: \PP\circ\PP\rto \PP$ and a unit
$u_\PP: I\rto \PP$ satisfying
\begin{align*}
\mu_\PP(\PP\circ\mu_\PP)=&\mu_\PP(\mu_\PP\circ\PP) \\
\mu_\PP(\PP\circ u_\PP)=&\mu_\PP( u_\PP\circ \PP)=\PP.
\end{align*}

In other terms an operad $\PP$ is a collection of $S_n$-modules
$(\PP(n))_{n\geq 0}$
together with an element $1_1\in\PP(1)$ and compositions
$$\gamma:\PP(k)\otimes\PP(l_1)\otimes\ldots\otimes\PP(l_k)\rto
\PP(l_1+\ldots +l_k)$$
satisfying associativity, unitary conditions and
equivariance conditions reflecting
the action of $S_k\times S_n$ on $\PP^{\hat\otimes k}(n)$ that is:
\begin{align*}
\gamma(p\cdot\sigma,p_1,\ldots,p_k)=&\gamma(p,p_{\sigma^{-1}(1)},
\ldots,p_{\sigma^{-1}(k)})\cdot\sigma(l_1,\ldots,l_k),\\
\gamma(p,p_1\cdot\tau_1,\ldots,p_k\cdot\tau_k)=&
\gamma(p,p_1,\ldots,p_k)\cdot(\tau_1\oplus\ldots\oplus\tau_k).
\end{align*}
Most of the time $\gamma(p,p_1,\ldots,p_k)$ is written $p(p_1,\ldots,p_k)$.
When $p_j=1_1$ for all $j$ except $i$ the latter composition is written
$p\circ_i p_i$.

 With this notation the associativity, unitarity and
equivariance write, for $p\in \PP(n), q\in\PP(m),r\in \PP(l),
\sigma\in S_n, \tau\in S_m$:

\begin{align}
\label{Assoc1-operad}
(p\circ_i q)\circ_{j+i-1} r=& p\circ_i (q\circ_j r), \\
\label{Assoc2-operad}
(p\circ_i q)\circ_{j+m-1} r=& (p\circ_j r)\circ_i q,\quad i<j, \\
\label{Unite-operad}
p\circ_i 1_1=&p=1_1\circ_1 p, \\
\label{Equivariance-operad}
(p\cdot\sigma)\circ_i (q\cdot\tau)=& (p\circ_{\sigma(i)} q)\cdot
  (\sigma\circ_i \tau),
\end{align}
where $\sigma\circ_i \tau$ is the permutation of $S_{n+m-1}$ obtained
by replacing $\tau$ for $\sigma(i)$. For instance
$$(3,4,2,5,1)\circ_2 (a,b,c)=(3,a+3,b+3,c+3,2,7,1).$$

%%%%%%%%%%%%%% pre-unital and unital operad

\subsubsection{Pre-unital and connected operads}
A \sl pre-unital \rm operad is an injective species $\II$ together with an operad structure such that the product $\mu_{\II}:\II\circ \II\rightarrow \II$ is a morphism of $\I$-modules.

An operad $\PP$ is said to be {\sl connected} if
$\PP(0)=\kk$ (unital in the terminology of May \cite{May72}).
With $\PP$ connected, let $1_0$ be
the image of $1 \in\kk=\PP(0)$. The operad $\PP$ is endowed with
degeneracy maps $\partial_i$ for $1\leq i\leq n$ defined by
$$\begin{array}{cccc}
\partial_i: &\PP(n)&\rto& \PP(n-1) \\
&\mu &\mapsto& \mu\circ_i 1_0.
\end{array}.$$
The composition of degeneracies induces restriction maps: for
any subset $S\subset [n]$ of size $k$ one defines
$$|_S: \PP(n)\rto\PP(k)$$
by
$$\mu|_{S}=\begin{cases}\mu, & {\rm if}\ S=[n], \\
\partial_{t_1}\partial_{t_2}\ldots \partial_{t_{n-k}}(\mu),& {\rm otherwise},\end{cases}$$ where
$\{t_1<t_2<\ldots <t_{n-k}\}=[n]\setminus S$. Hence
$$\mu|_{S}=\mu(x_1,\ldots,x_n)\  {\rm where}\
  x_i=\begin{cases} 1_1,&
  \rm{if}\ i\in S, \\ 1_0, & \rm{otherwise}. \end{cases}$$

%%%%%%%%%%%%%%%%%% the associative operad

\subsubsection{The associative operad} \label{S-associativeoperad}
The operad $\As$ is defined by $\As(n)=\kk[S_n]$ for all $n\geq 0$.
The composition $\sigma\circ_i\tau$ is the one given in definition \ref{D-operad}.
This operad is connected.
The degeneracy maps were defined in paragraph \ref{S-symmetricgroups}: 
$$\partial_i(\sigma)=\sigma\circ_i 1_0=
\st(\sigma_1,\ldots,\sigma_{i-1},\sigma_{i+1},\ldots,\sigma_n).$$
For $S=\{s_1<\ldots<s_k\}$ the permutation  
$\sigma|_{S}$ is the standardisation of 
$(\sigma(s_1),\ldots,\sigma(s_k))$. For instance
$$(3,2,6,1,8,7,5,4)|_{\{1,4,6,7\}}=\st(3,1,7,5)=(2,1,4,3).$$
Furthermore, for any $\sigma\in S_n$
$$\sigma|_{\emptyset}=1_0.$$

%%%%%%%%%%%%%%%%%%%%%%% toute operade unitaire est une operade pre-unitale

\subsubsection{Proposition}\it Any connected operad is a pre-unital operad.\rm

\medskip

\noindent{\sl Proof.}  By definition $\partial_i(\mu)=\mu\circ_i 1_0$. The formula
(\ref{Assoc2-operad}) implies equation (\ref{E-degeneracy1}) and the formula
 (\ref{Equivariance-operad}) implies equation (\ref{E-degeneracy2}).
This proves that $\PP$ is a
$\I$-module.

Let $\mu\in\PP(n), \nu_k\in \PP(l_k)$ and $i\in [R_{j-1}+1,R_{j-1}+l_j]$ with 
$R_{j-1}=l_1+\ldots+l_{j-1}$. Iterating formula 
(\ref{Assoc1-operad}) gives
\begin{multline*}
\partial_i(\mu_{\PP}(\mu\otimes\nu_1\otimes\ldots\otimes\nu_n))=\mu(\nu_1,\ldots,\nu_{j-1},
\partial_{i-R_{j-1}}(\nu_j),\ldots,\nu_n) \\
=\mu_{\PP}(\partial_i(\mu\otimes\nu_1\otimes\ldots\otimes\nu_n)).
\end{multline*}
Hence $\mu_{\PP}$ is a morphism of $\I$-modules.  \hfill$\Box$

\medskip

Given $n$ sets $S_i\subset [l_i]$ for $1\leq i\leq n$ the set
$S_1\star\ldots\star S_n$ is the subset of $[l_1+\ldots+l_n]$ of elements $R_{i-1}+\alpha$ for
$\alpha\in S_i$.
Iterating the last equation gives the relation
$$\mu(\nu_1,\ldots,\nu_n)|_{S_1\star\ldots\star S_n}=\mu(\nu_1|_{S_1},\ldots,\nu_n|_{S_n}),$$
which yields the following lemma

%%%%%%%%%%%%%%%%%%%%%% Corollaire: des formules

\subsubsection{Lemma}\label{L-partialonsets}\it Let $\PP$ be a connected operad.
Let $\mu\in\PP(n),\nu_i\in\PP(l_i)$ and $S_i\subset [l_i]$. For any set
$J=\{j_1<\ldots <j_{l}\}\subset [n]$ such that $S_i$ is empty for all $i\not\in J$, one has

\begin{equation*}
\mu(\nu_1,\ldots,\nu_n)|_{S_1\star\ldots\star S_n}=
(\prod_{i\not\in J} \nu_i|_{\emptyset}) \mu|_J(\nu_{j_1}|_{S_{j_1}},\ldots,\nu_{j_{l}}|_{S_{j_{l}}}).
\end{equation*}\rm

%%%%%%%%%%%%%  SUBSECTION %%%%%%%%%%%%%%%%%% TWISTED ALGEBREAS OVER AN OPERAD

\subsection{Twisted algebras over an operad}
  Let $V$ be a vector space. It can be considered as a $\S$-module concentrated in degree 0. 
In particular, the notation $\PP\circ V$ makes sense, and an algebra over an operad is a vector 
space together with a product map $\PP\circ V\rto V$ satisfying the usual monadic relations 
(see Ginzburg and Kapranov  \cite{GinKap94}  or replace $M$ by $V$ in the diagrams below).
The notion of twisted algebras over an operad generalizes
this definition to $\S$-modules. Notice that, since the category of vector spaces identifies 
with the full subcategory of the category of $\S$-modules which objects are the $\S$-modules concentrated 
in degree 0, any result on twisted algebraic structures holds automatically for classical 
algebraic structures. We will refer to this property as ``Restriction to $\Vect$''.

%%%%%%%% definition

\subsubsection{Definition} \label{D-twistedalgebra} Let $\PP$ be an operad. A
$\S$-module $M$ is a {\sl twisted $\PP$-algebra} if
$M$ is endowed with a product
$\mu_M: \PP\circ M\rto M$ such that the following diagrams commute:

$$\xymatrix{\PP\circ\PP\circ M \ar_{\mu_\PP\circ M}[d] \ar^>>>>>{\PP\circ\mu_M}[r]
& \PP\circ M \ar^{\mu_M}[d] \\
\PP\circ M\ar^{\mu_M}[r] & M} \qquad \qquad\xymatrix {I\circ M\ar_{ u_\PP\circ M}[d]
\ar^{=}[r] & M\\
\PP\circ M \ar_{\mu_M}[ur] &}$$

%%%%%%%%%%%%%%%%%%%% remark

\subsubsection{Remark} Twisted $\PP$-algebras are called left $\PP$-modules
in the categorical, monadic, terminology
(see e.g. \cite{Fr04}), since they are also algebras over the monad
$$\begin{array}{cccc}
\PPP:& \S-mod & \rto & \S-mod \\
& M & \mapsto & \PP\circ M
\end{array}
$$

Here, we prefer to stick to the more appealing terminology of twisted $\PP$-algebras. 
In case $\PP$ is the operad defining associative or Lie algebras,
the definition coincides with the notion of
twisted associative or twisted Lie algebras in the sense of Barratt in
\cite{Barratt78}, see the example below and Joyal's article \cite{Joyal86}.

If $M=M(0)$ is concentrated in
degree 0 then $M$ is a $\PP$-algebra in the usual sense,   according to the Restriction to 
$\Vect$ principle. For instance
if $\PP=\As$, $\As\circ M=T(M)$ is the free associative algebra over $M$.
If $M$ is a $\S$-module such that the action of $S_n$ is the trivial action or the signature
action, then $M$ is a graded $\PP$-algebra in the usual sense (when the action is trivial no
signs are considered).

%%%%%%%%%%%%%%%%%%%%%%%%%% example: twisted Lie algebras

\subsubsection{Example: Twisted Lie algebras} Let $\Lie$ be the Lie operad. As an operad it is generated
by $\mu\in\Lie(2)$ satisfying the following relations:
\begin{align*}
&\mu\cdot(2,1)=-\mu \\
&\mu\circ_2\mu\cdot( (1,2,3)+(2,3,1)+(3,1,2))=0
\end{align*}
A {\sl twisted Lie algebra} is a twisted algebra over the operad $\Lie$, that is a
$\S$-module $M$ endowed with a multiplication, the bracket,
defined by $[a,b]=\mu(a,b)$ satisfying the following relations: let
$a\in M(p), b\in M(q), c\in M(r)$,
\begin{align*}
&[b,a]\cdot\zeta_{p,q}=-[a,b] \\
&[a,[b,c]]+[c,[a,b]]\cdot \zeta_{p+q,r}+[b,[c,a]]\cdot \zeta_{p,q+r}=0,
\end{align*}
where $\zeta_{p,q}$ was defined in relation (\ref{E-symmetryiso}).
This is exactly definition 4 in \cite{Barratt78}.
These relations are a direct  consequence of the computation of $(\mu\cdot(2,1))(a,b)$ and
$((\mu\circ_2\mu)\cdot (1,2,3)+(2,3,1)+(3,1,2))(a,b,c)$ using the symmetry isomorphism $\tau$.
For instance
$$((\mu\circ_2\mu)\cdot (2,3,1))(a,b,c)=[c,[a,b]]\cdot(2,3,1)(p,q,r)=[c,[a,b]]\cdot \zeta_{p+q,r}.$$

\subsubsection{Unital $\PP$-algebras}\label{D-unitalalgebra} Assume $\PP$ is a connected operad, then the $\S$-module $\un$ is a
twisted $\PP$-algebra for the product
$$\mu(1_0,\ldots,1_0)=\mu|_{\emptyset}1_0.$$
By definition a {\sl unital twisted $\PP$-algebra} is a twisted $\PP$-algebra $M$ together with a morphism of
twisted $\PP$-algebras called the unit morphism
$$\eta_M: \un\rto M$$
satisfying
\begin{equation}\label{R-unitalalgebra}
\mu(a_1,\ldots,a_n)=\mu|_{\{i_1,\ldots,i_r\}}(a_{i_1},\ldots,a_{i_r})
\end{equation}
for all $\mu\in\PP(n), a_i\in M$ as soon as $a_j=\eta_M(1_0)$ for all $j\not\in \{i_1,\ldots,i_r\}$.
A {\sl morphism of  unital twisted $\PP$-algebras} is a morphism of twisted $\PP$-algebras commuting with
the unit morphism.

%%%%%%%%%%%%%%%%%%%%%%%%%%%%%% the free twisted P-alegbra on one generator.

\subsubsection{Free twisted $\PP$-algebras}\label{S-free} Since $I$ is the unit
for the plethysm, the $\S$-module $\PP=\PP\circ I$ is the free twisted
$\PP$-algebra generated by $I$.
More generally, $\PP\circ M$ is the free twisted $\PP$-algebra
generated by a $\S$-module $M$. If $\PP$ is connected then $\PP\circ
M$ is a unital twisted  $\PP$-algebra since

$$\eta_{\PP\circ M}(0): \kk\rto (\PP\circ M)(0)=\PP(0)\oplus_{k\geq 1}\PP(k)\otimes_{S_k} M(0)^{\otimes k}.$$
is given by the isomorphism $\PP(0)=\kk$.

%%%%%%%%%%%%%%%%%%%%%%  SECTION %%%%%%%%%%%%  TWISTED ALGEBRAS OVER HOPF OPERAD

\section{twisted Hopf algebras over a Hopf operad}

%%%%%%%%%%%%%%%%%%%%% SUBSECTION %%%%%%%%%%%%%%% HOPF OPERADS

\subsection{Hopf operads}

Let $\Coalg$ be the category of coassociative counital coalgebras, that is vector spaces $V$ endowed with
a coassociative coproduct $\delta: V\rto V\otimes V$ and
a linear map, the counit,
$\epsilon: V\rto \kk$ satisfying
$$(\epsilon\otimes V)\delta=(V\otimes\epsilon)\delta=V.$$

%%%%%%%%%%%%%%%%% definition: hopf operads

\subsubsection{Definition} \label{D-hopfoperad} A {\sl Hopf operad $\PP$} is an operad
in the category of coalgebras: $\mu_\PP$ and $ u_\PP$ are morphisms of coalgebras.
More precisely, for each $n$ there exists a coassociative coproduct
$$\delta(n):\PP(n)\rto \PP(n)\otimes\PP(n)=(\PP\otimes\PP)(n)$$
which is $S_n$-equivariant and such that the following two diagrams commute,
$$\xymatrix{ \PP\circ \PP \ar_{\delta\circ\delta}[d]\ar^{\mu_P}[rr] & & \PP \ar_{\delta}[d] \\
(\PP\otimes\PP)\circ(\PP\otimes\PP)\ar^{T_1}[r]& (\PP\circ\PP)\otimes(\PP\circ\PP)
\ar^>>>>{\mu_\PP\otimes\mu_\PP}[r] &\PP\otimes\PP}$$
$$\xymatrix{I \ar_{\delta_I}[d]\ar^{ u_\PP}[r]& \PP \ar_{\delta}[d] \\
I\otimes I\ar^{ u_P\otimes u_\PP}[r]& \PP\otimes\PP}
$$
and for each $n$ there exists
$\epsilon(n):\PP(n)\rto \kk$ (abbreviated to $\epsilon$ when no confusion can arise) such that
\begin{align*}
\epsilon(\mu(\nu_1,\ldots,\nu_n))=&\epsilon(\mu)\epsilon(\nu_1)\ldots\epsilon(\nu_n)\\
\epsilon(1_1)=&1_{\kk}.
\end{align*}
A {\sl connected Hopf operad} is a Hopf structure on a connected operad such that the counit
$\epsilon$ is given by
$$\epsilon(\mu)=\mu|_\emptyset 1_{\kk}.$$

%%%%%%%%%%%%%%%%%%%%%%%%%%%Example

\subsubsection{Example} The operad $\Com$ is defined by $\Com(n)=\kk$ with the trivial action of the symmetric group.
The composition map is given by $e_n\circ_i e_m=e_{n+m}$ where $e_n$ states for $1_{\kk}\in\Com(n)$. It is a connected
Hopf operad for the coproduct $\delta(e_n)=e_n\otimes e_n$ and the counit $\epsilon(e_n)=1_{\kk}$.

The operad $\As$ is a connected Hopf operad for the coproduct
$\delta(\sigma)=\sigma\otimes\sigma$ and $\epsilon(\sigma)=1_{\kk}$ for all $\sigma\in S_n$.

%%%%%%%%%%%%%%%%%%%%%%%%%%%% theorem : si P est de Hopf, les P-algebres ont un produit tensoriel

\subsubsection{Theorem}\label{T-hopfhatotimes}\it Let $\PP$ be a Hopf operad and $V,W$ be twisted $\PP$-algebras. The
$\S$-module $V\hat\otimes W$ is a twisted $\PP$-algebra for the
following product: $\mu_{V\hat\otimes W}$ is the composite of
$$\xymatrix{
\PP\circ (V\hat\otimes W)\ar^>>>>{\delta\circ(V\hat\otimes W)}[r]& (\PP\otimes\PP)\circ
(V\hat\otimes W) \ar^{T_2}[r] & (\PP\circ V)\hat\otimes (\PP\circ W)
\ar^>>>>{\mu_V\hat\otimes\mu_W}[r]& V\hat\otimes W}.$$
If the coalgebra structure on $\PP$ is cocommutative, then the symmetry isomorphism $\tau$ is a morphism
of $\PP$-algebras.
\rm

\medskip

\noindent{\sl Proof.} Let $\delta:\PP\rto \PP\otimes\PP$ be the
Hopf coproduct of $\PP$ and  $M:=V\hat\otimes W$. One has to prove that
the following diagrams commute:

$$\xymatrix{\PP\circ\PP\circ M \ar_{\mu_\PP\circ M}[d] \ar^>>>>>{\PP\circ\mu_M}[r]
& \PP\circ M \ar^{\mu_M}[d] \\
\PP\circ M\ar^{\mu_M}[r] & M} \qquad \qquad\xymatrix {I\circ M\ar_{\eta_\PP\circ M}[d]
\ar^{=}[r] & M\\
\PP\circ M \ar_{\mu_M}[ur] &}$$

For the first diagram, one has

\begin{align*}
\mu_M(\mu_\PP\circ M)=& (\mu_V\hat\otimes
\mu_W)T_2(\delta\circ(V\hat\otimes W))(\mu_\PP\circ(V\hat\otimes W)) \\
=& (\mu_V\hat\otimes
\mu_W)T_2((\delta\mu_\PP)\circ(V\hat\otimes W)).
\end{align*}
The following relations hold
\begin{align*}
(\delta\mu_\PP)\circ (V\hat\otimes
W)=&(\mu_\PP\otimes\mu_\PP)T_1(\delta\circ\delta)\circ (V\hat\otimes
W) \\
T_2((\mu_\PP\otimes \mu_\PP)\circ V\hat\otimes W)&=((\mu_\PP\circ
V)\hat\otimes (\mu_\PP\circ W))T_2 \\
(\mu_V\hat\otimes\mu_W)((\mu_\PP\circ V)\hat\otimes(\mu_\PP\circ
W))&=(\mu_V\hat\otimes \mu_W)((\PP\circ \mu_V)\hat\otimes(\PP\circ \mu_W))\\
\end{align*}
The first equality comes from definition \ref{D-hopfoperad}. The
second one comes from the naturality of $T_2$. The last one comes 
from the definition of twisted $\PP$-algebras applied on $V$ and $W$
(see \ref{D-twistedalgebra}).
Combining these equalities, one gets
\begin{multline*}
\mu_M(\mu_\PP\circ M)=\\
=(\mu_V\hat\otimes \mu_W)((\PP\circ \mu_V)\hat\otimes(\PP\circ
\mu_W))T_2(T_1\circ M)((\delta\circ\delta)\circ M)
\end{multline*}
The relation (\ref{R_T1_T2}) implies
\begin{align*}
T_2(T_1\circ M)((\delta\circ\delta)\circ M)=&T_2((\PP\otimes\PP)\circ
T_2)((\delta\circ\delta)\circ M) \\
\end{align*}
The naturality of $T_2$ implies
\begin{align*}
((\PP\circ \mu_V)\hat\otimes(\PP\circ
\mu_W))T_2=&T_2((\PP\otimes\PP)\circ(\mu_V\hat\otimes\mu_W))
\end{align*}
As a consequence one has
\begin{align*}
\mu_M(\mu_\PP\circ M)=&
(\mu_V\hat\otimes
\mu_W)T_2((\PP\otimes\PP)\circ(\mu_V\hat\otimes\mu_W))
((\PP\otimes\PP)\circ T_2)((\delta\circ\delta)\circ M)\\
=&(\mu_V\hat\otimes
\mu_W)T_2(\delta\circ(\mu_V\hat\otimes\mu_W))
(\PP\circ T_2)((\PP\circ\delta)\circ M)\\
=&(\mu_V\hat\otimes
\mu_W)T_2(\delta\circ M)(\PP\circ \mu_M)=\mu_M(\PP\circ\mu_M).
\end{align*}
It is easy to check that $\mu_M$ makes the second diagram
commute. 

In case $\delta$ is cocommutative, using diagrams one sees easily that
$\tau: V\hat\otimes W\rto W\hat\otimes V$ is a morphism of $\PP$-algebras.
\hfill $\Box$

\medskip

%%%%%%%%%%%%%% Rem: lien avec Moerdijk
One could deduce the theorem from Moerdijk 
\cite[Proposition 1.4]{Moerdijk02}, since
a part of the proof above implies that the monad
$$\begin{array}{cccc}
\PPP:& \S-mod & \rto & \S-mod \\
& M & \mapsto & \PP\circ M
\end{array}
$$
is a {\sl Hopf monad}. However the proof above is still valid if we
require the operad to be an operad in the category of coassociative
coalgebras (non necesseraly counital). The counit is necessary for
the category of twisted $\PP$-algebras to be a tensor category: the unit for $\hat\otimes$
is given by the $\PP$-algebra $\un$ endowed with the product
$$\mu_{\un} :(\PP\circ \un)(\emptyset)=\oplus_{l\geq 0} \PP(l)/S_l \rto \un(\emptyset)=\kk$$
given by $\mu_{\un}(\nu)=\epsilon(n)(\nu)$ for $\nu\in\PP(n)$.

When the operad is a Hopf connected operad the $\PP$-algebra structure on $\un$ 
coincides with the one given in \ref{D-unitalalgebra}. This implies the following

\subsubsection{Theorem}\it Let $\PP$ be a connected Hopf operad. Then the category of unital
twisted $\PP$-algebras is a tensor category, symmetric if $\PP$ is cocomutative. \rm

%%%%%%%%%%%%%% SUBSECTION %%%%%%%%%% TWISTED ALGEBRAS OVER HOPF OPERADS

\subsection{Twisted Hopf algebras over a connected Hopf operad}

%%%%%%%%%%%%%%%%%%%%% definition

\subsubsection{Definition} Let $\PP$ be a connected Hopf operad. A
unital twisted $\PP$-algebra $M$ is a
{\sl twisted Hopf $\PP$-algebra} if $M$ is endowed with a coassociative counital
coproduct
$$\Delta: M\rto M\hat\otimes M, \qquad  \epsilon_M: M\rto \un$$
where $\Delta$ and $\epsilon_M$ are morphisms of unital $\PP$-algebras.

%%%%%%%%%%%%%%%%%%%%%%% remarque: as

\subsubsection{Remark} In case $\PP=\As$ one gets twisted bialgebras as defined in 
\cite{Stover93,PatReu04}.

%%%%%%%%%%%%%%%%%%%% theorem: P unitale P est twisted Hopf.

\subsubsection{Theorem}\label{T-hopfhopf}{\it If $\PP$ is a connected Hopf operad
then $\PP$ is a twisted Hopf $\PP$-algebra for
the coproduct
\begin{equation}\label{E-DeltaOperade}
\Delta(\mu)=\sum_{{(1),(2)}\atop{S\sqcup T=[n]}}
\mu_{(1)}|_S\otimes \mu_{(2)}|_T \cdot \sigma(S,T)^{-1}
\end{equation}
where the Hopf structure on $\PP$ is given by
$\delta(\mu)=\sum_{(1),(2)} \mu_{(1)}\otimes \mu_{(2)}$
and for $S=\{s_1<\ldots<s_k\}$ and $T=\{t_1<\ldots<t_{n-k}\}$ the permutation
$\sigma(S,T)$ is $(s_1,\ldots,s_k,t_1,\ldots,t_{n-k})$.
The coproduct $\Delta$ is cocommutative if $\PP$ is a cocommutative Hopf operad.
}

\medskip

{\sl Proof.} Consider the map of $\S$-modules induced by the embeddings $I\oplus I\subset \PP(1)\oplus \PP(1)=\PP (1)\otimes \PP(0)\oplus \PP(0)\otimes\PP(1)\subset \PP\hat\otimes \PP$:

$$\begin{array}{cccc}
\phi: &I&\rto& \PP\hat\otimes \PP \\
& 1_1 & \mapsto & (1_1\otimes 1_0)+(1_0\otimes 1_1).\end{array}$$

Since $\PP$ is the free twisted $\PP$-algebra on $I$ and since
$\PP\hat\otimes\PP$ is endowed with a twisted $\PP$-algebra structure
thanks to theorem \ref{T-hopfhatotimes}, there is a unique morphism of twisted
$\PP$-algebras 
$$\Phi:\PP\circ I=\PP \rto \PP\hat\otimes\PP$$
extending $\phi$. Indeed
$$\Phi=\mu_{\PP\hat\otimes\PP}(\PP\circ\phi)=(\mu_\PP\hat\otimes\mu_\PP)T_2(\delta\circ\phi)$$
As a consequence
\begin{align*}
\Phi(\mu)&=(\mu_\PP\hat\otimes\mu_\PP)T_2(\sum_{(1),(2)}\mu_{(1)}\otimes\mu_{(2)}\otimes
(1_0\otimes
1_1+1_1\otimes 1_0,\ldots,1_0\otimes
1_1+1_1\otimes 1_0)\\
&=(\mu_\PP\hat\otimes\mu_\PP)(\sum_{{(1),(2)}\atop {S\sqcup
  T=[n]}}\mu_{(1)}\otimes (x_1\otimes\ldots x_n)\otimes\mu_{(2)}\otimes
(y_1\otimes\ldots y_n) \cdot\sigma(S,T)^{-1}) \\
&=\sum_{{(1),(2)}\atop{S\sqcup T=[n]}}
\mu_{(1)}|_S\otimes \mu_{(2)}|_T \cdot \sigma(S,T)^{-1}
\end{align*}
where $$\begin{cases} x_i=1_1\ {\rm and}\  y_i=1_0,\ &{\rm if}\ i\in S,
  \\  x_i=1_0 \ {\rm and}\
  y_i=1_1,\ &\ {\rm if}\ i\in T,\end{cases}$$
and $\sigma_{S,T}^{-1}$ is the shuffle coming from $T_2$.
The unit $\eta_{\PP}$ and counit $\epsilon_{\PP}$ 
are given by the isomorphism between $\kk$ and $\PP(0)$. Indeed
$$(\epsilon_{\PP}\otimes \Id)\Delta(\mu)=\sum_{(1),(2)}(\mu_{(1)}|_{\emptyset})\mu_{(2)}=\mu$$
because $\mu_{(1)}|_\emptyset =\epsilon(\mu_{(1)})$ and $\epsilon$ is the counit for $\delta$.

The coasociativity of $\Delta$ follows from the coassociativity of $\phi$ and the unicity of
the extension of $\phi$ as a morphism of $\PP$-algebras. The cocommutativity is clear in case $\PP$ is
a cocommutative Hopf operad. \hfill$\Box$

%%%%%%%%%%%%%%%%%%%%% Example: as

\subsubsection{Example} In case $\PP=\As$ 
one has
$$\Delta(\sigma)=\sum_{S\sqcup T=[n]}
\sigma|_{S}\otimes \sigma|_{T} \cdot \sigma(S,T)^{-1}.$$

This is the twisted Hopf algebra structure on the direct sum of the symmetric group 
algebras described by Patras and Reutenauer in \cite{PatReu04}.

%%%%%%%%%%%%%%%%%%%% corollary

\subsubsection{Corollary}\it Let $\PP$ be a connected Hopf operad. Any free 
twisted $\PP$-algebra has a canonical twisted Hopf $\PP$-algebra structure,
cocommutative whenever $\PP$ is. \rm

\medskip

\noindent{\sl Proof.} The map $\PP\circ M\rto (\PP\circ M)\hat\otimes (\PP\circ M)$ is given
by the composite
$$\PP\circ M \rto (\PP\hat\otimes \PP)\circ M =(\PP\circ M)\hat\otimes (\PP\circ M).$$
thanks to relation (\ref{R-iso}). The unit $\eta_{\PP\circ M}$ and counit $\epsilon_{\PP\circ M}$ 
are induced by  the one on $\PP$: since $\un\circ M=\un$,
one has $\eta_{\PP\circ M}:=\eta_{\PP}\circ M: \un=\un\circ M\rto \PP\circ M$ and
$\epsilon_{\PP\circ M}:=\epsilon_{\PP}
\circ M: \PP\circ M\rto \un\circ M=\un$. \hfill$\Box$

%%%%%%%%%%%%%%%%%%Primitive elements

\subsection{Primitive elements}\label{S-Prim}

%%%%%%%%%%%%Definition

\subsubsection{Definition}\label{D-primitive} Let $\PP$ be a connected Hopf operad and  $M$ be a twisted Hopf
$\PP$-algebra. Then $M$ has a unit $\eta_M:\un\rto M$ and  a counit $\epsilon_M:M\rto\un$ such
that $\epsilon_M\eta_M=\un$. Let
$1_M$ be the image of $1_0$ by $\eta_M$.
Let $\overline M:=\Ker\epsilon_M$. One has
$$M=\un\oplus\overline M$$
and for all $x\in\overline M$ one has
$$\Delta_M(x)=1_M\otimes x+x\otimes 1_M+\sum_{i,j} x_i\otimes x_j$$
with $x_i,x_j\in \overline M$.
The space $\Prim(M)$
of {\sl primitive elements} of $M$  is
$$\Prim(M)=\{x\in\overline M| \Delta(x)=1_M\otimes x+x\otimes 1_M\}.$$
Since $\Delta_M$ is a morphism of $\S$-modules, $\Prim(M)$ is a sub $\S$-module of $M$.
In the sequel, $\overline\Delta$ denotes the projection of $\Delta$ onto 
$\overline M\hat\otimes \overline M$.
The space of primitive elements is then 
$\Prim(M)=\ker(\overline\Delta:\overline M\rto \overline M\hat\otimes \overline M$).
Note that if $V$ is a $\S$-module then $\PP\circ V$ is a twisted Hopf $\PP$-algebra and
\begin{equation}\label{overline}
\overline{\PP\circ V}=\overline{\PP}\circ V
\end{equation}

\medskip

Since $\overline\Delta$ is coassociative, one can define consistently 
$\overline\Delta^{[n]}$ as a map from $\overline M$ to $\overline M^{\hat\otimes n}$ 
(e.g. $\overline\Delta^{[3]}:=\overline\Delta\hat\otimes \overline M =
\overline M\hat\otimes\overline\Delta$). The twisted Hopf $\PP$-algebra $M$
is said to be {\sl connected} if for any $x\in \overline M$ there exists $n$ such that
$$\overline \Delta^{[n]}(x)=0.$$ For instance, $\PP$ is a connected twisted Hopf
$\PP$-algebra.

%%%%%%%%%%%%%%%%Theoreme fondamental

\subsubsection{Theorem}\label{FundamentalPrim}{\it  Let $\PP$ be a connected Hopf operad.
The space of primitive elements of the twisted Hopf $\PP$-algebra $\PP$
is a sub-operad of $\PP$.}

\medskip

{\sl Proof.} Notice that $\mu\in\PP(n)$ is primitive if and only if
\begin{equation*}
\sum_{S\sqcup T=[n], S,T\not=\emptyset} \mu_{(1)}|_S\otimes
\mu_{(2)}|_T\cdot \sigma(S,T)^{-1}=0,
\end{equation*}
that is if and only if $\Delta_{S,T}(\mu):=\mu_{(1)}|_S\otimes
\mu_{(2)}|_T=0$.

As pointed out in the definition \ref{D-primitive}, the space  $Q=\Prim(\PP)$ is a
sub-$\S$-module of $\PP$; moreover $1_1\in Q(1)$.

Assume $p\in Q(n)$, $q_i\in Q(l_i), l_i>0$ for all $1\leq i\leq n$ and
$S\sqcup T=[l_1+\ldots+l_n], S,T\not=\emptyset$.
Let us write $S=S_1\star\ldots\star S_n$ and $T=T_1\star\ldots\star T_n$ 
with $S_i\sqcup T_i=[l_i]$.
Since $\Delta:\PP\rto\PP\hat\otimes\PP$ is a morphism of $\PP$-twisted
algebras one has

\begin{multline*}
\Delta_{S,T}(p(q_1,\ldots,q_n))= \\
\sum_{(a),(b),(a_1,\ldots,a_n),(b_1,\ldots,b_n)}
  p_{(a)}(q_{1(a_1)}{|_{S_1}},\ldots ,q_{n(a_n)}{|_{S_n}})\otimes
  p_{(b)}(q_{1(b_1)}{|_{T_1}},\ldots ,q_{n(b_n)}{|_{T_n}}).
\end{multline*}

Since the $q_i$ are primitive elements, the displayed
quantity is equal to $0$ if there exists $i\leq n$ with
$S_i\not= \emptyset $ and $T_i\not= \emptyset$. So, let us assume
that, for any $i$, $S_i=\emptyset $ or
 $T_i=\emptyset$, and let us set: $\{
 i_1,...,i_k\}=\{i,S_i\not=\emptyset\}$ and $\{
 j_1,...,j_{n-k}\}=[n]-\{i_1,...,i_k\}=
\{j,T_j\not=\emptyset\}$.
We get by lemma \ref{L-partialonsets}
\begin{multline*}
\Delta_{S,T}(p(q_1,\ldots,q_n))=\\
\sum\limits_{(1),(2)}p_{(1)}|_{\{i_1,...,i_k\}}(q_{i_1},...,q_{i_k})\otimes p_{(2)}|_{\{j_1,...,j_{n-k}\}}(q_{j_1},...,q_{j_{n-k}})
\end{multline*}
which is equal to zero since $p$ is primitive, and $S,T\not=\emptyset$
implies that $k,n-k\not= 0$. \hfill $\Box$

\medskip

For example, if $\PP =\As$, it follows from the theorem and from
\cite[Prop.17]{PatReu04} that the operad of
 primitive elements of $\As$, viewed as a twisted Hopf algebra, is the
 Lie operad.
This is not a surprising result in view of the classical Lie theory
and structure theorems for Hopf algebras such as
the Cartier-Milnor-Moore theorem \cite{MilMOO65,Patras94}, however it
shows that many properties of classical Hopf algebraic
structures are encapsulated in the twisted Hopf algebra structure of the corresponding Hopf operad.

The next sections are devoted to the systematical study of these properties.

Before turning out to this systematical study, let us consider the example of the magmatic operads 
considered by Holtkamp in \cite{Hol05}.

\subsection{Magmatic operads} In \cite{Hol05}, Holtkamp considers free operads \hfill\break
$\Mag_N$
and $\Mag_\omega$ which are generated by an operation in $k$ variables
$\vee^k$ for each $2\leq k\leq N$ or for
each $2\leq k$. These operads can be turned out to connected operadsby setting
$\vee^k\circ_i 1_0=\vee^{k-1}$ and more generally $\vee^k|_S=\vee^{|S|}$, where $\vee^1=1_1$ and $\vee^0=1_0$. Moreover these operads
are connected Hopf operad with $\delta(\vee^k)=\vee^k\otimes\vee^k$ and $\epsilon(\vee^k)=1_{\kk}$.
As a consequence, we recover from the previous theorem some results of his paper. For instance, $\Prim(\Mag_N)$ is a suboperad
of $\Mag_N$ and the same property holds for $\Mag_{\omega}$.

%%%%%%%%%%%%%%%%%%%%%%%%%%%%%%%Lois de reciprocite

\section{Reciprocity laws}

\subsection{Hopf algebras over an operad}
In the present section, we investigate the relations between the structure of a connected 
Hopf operad and the
structure of Hopf algebras (twisted and classical) over this operad.

In this section, and the following one, $\PP$ is a connected Hopf operad, and $Q$ is the operad of
primitive elements of $\PP$.

%%%%%%%%%%%%%%%%% Theorem: Prim(H) est une Q-algebre.

\subsubsection{Theorem}\label{Prim} \it Let $H$ be a twisted Hopf $\PP$-algebra.
Then, the $\S$-module $\Prim(H)$ is a $Q$-algebra.\rm

\medskip

\noindent{\sl Proof.} Since any element in $\Prim(H)$ satisfies $\Delta_H(h)=1_H\otimes h+h\otimes 1_H$
the same proof as in theorem \ref{FundamentalPrim} holds using relation (\ref{R-unitalalgebra}) instead of lemma
\ref{L-partialonsets}. \hfill$\Box$

\medskip

As a direct consequence, due to the Restriction to $\Vect$ principle, 
the set of primitive elements of any Hopf algebra over $\PP$
is a $\Prim(\PP)$-algebra.

%%%%%%%%%%%%%%%%%%%%%%%%%%%%% Lie monomial

\subsubsection{Generalized Lie monomials}

Recall that, in the classical theory of free Lie algebras, a \it Lie
polynomial \rm in a free associative algebra (or tensor algebra)
$T(V)$ over a vector space $V$ is an arbitrary element in the free Lie
algebra over $V$, where the latter is viewed as a sub
Lie algebra of $T(V)$.

Due to a slightly misleading convention, a Lie monomial is a non
commutative monomial of Lie polynomials. Lie monomials and
Lie polynomials are one of the basic tools in the study of free Lie
algebras, and, actually, most of the properties of free Lie algebras
can be deduced from the behaviour of Lie monomials. We refer to the
book by Reutenauer \cite{Reu93}, where this point of view is
developped
in a systematic way.
More generally, Lie monomials can be defined in an arbitrary enveloping algebra $U(L)$, 
as non commutative monomials on $L$ (where, as usual, $U(L)$ is viewed as a 
quotient of $T(L)$).

From the combinatorial point of view, the fundamental property of Lie
monomials is their behaviour with respect to the coproduct in the
tensor
algebra and, more generally, in an arbitrary graded connected
cocommutative Hopf algebra: besides \cite{Reu93}, we refer to
\cite{Patras94},
and the computations therein, for more infomations on the subject.

The purpose of the next paragraphs is to extend this property to arbitrary Hopf operads.

%%%%%%%%%%Q-monomials

\subsubsection{Definition}
Let $H$ be a twisted Hopf $\PP$-algebra and $\Prim(H)$ the
$Q$-algebra of primitive elements in $H$. By analogy with the case of
Lie
monomials, we call {\sl $Q$-monomials} the elements of the free $\PP$-algebra over $\Prim(H)$.

\subsubsection{Theorem}\label{reciprocity}\it
We have the identity for $Q$-monomials:
\begin{multline}\label{PrimThm}
\forall \mu\in \PP (n), ~\forall h_1,...,h_n\in \Prim(H)\\
\Delta (\mu (h_1,...,h_n))=\sum\limits_{S,T}\Delta_{S,T}(\mu )(h_{i_1},...,h_{i_k}\otimes h_{j_1},...,h_{j_{n-k}})
\end{multline}
where $S=\{i_1,...,i_k\}$ and $T=\{j_1,...,j_{n-k}\}$ run over the
partitions of $[n]$, and where we write $\Delta_{S,T}(\mu )$ for the
component
of $\Delta (\mu )$ in $\PP(S)\otimes\PP(T)\subset (\PP\hat\otimes\PP )(n)$.\rm

\medskip

Once again, the theorem follows by adapting the proof of \ref{FundamentalPrim}.

Due to the Reduction to $\Vect$ principle, the theorem includes for example, as a particular case, 
the computation of coproducts of Lie monomials in an enveloping algebra. The
interested
reader may check that one recovers, in that case, the unshuffling coproduct
formula familiar in Lie theory \cite{Reu93}.

\medskip

%%%%%%%%%%%%%%%%%%%%%%%%%%%%%%Algebres libres

\subsection{Free algebras}
In the present section, we assume that $\kk$ is a field of characteristic zero.

%%%%%%%%%%%%%%%%Theorem Prim(P\circ V)=\Prim(P)\circ V

\subsubsection{Theorem}\it The $Q$-algebra of primitive elements of the
free twisted Hopf $\PP$-algebra $\PP\circ V$ over a  
$\S$-module $V$ is canonically isomorphic to the free twisted $Q$-algebra $Q\circ V$.\rm

\medskip

In particular, due to the Reduction to $\Vect$ principle, it follows that the primitive 
elements of the free
$\PP$-algebra over a vector space, viewed as a twisted Hopf $\PP$-algebra,
identify with the elements of the free $Q$-algebra over $V$. The
result generalizes to arbitrary Hopf operads the fundamental property
of associative algebras: the primitive part of a tensor algebra (i.e. of
a free associative algebra, naturally provided with a
cocommutative Hopf algebra structure) is a free Lie algebra \cite{Reu93}.

\medskip

\noindent{\sl Proof}. By definition of $Q$, the following sequence of $\S$-modules is left exact:
$$Q\mapright i \overline\PP\mapright{\overline\Delta}\overline{\PP}{\hat\otimes}\overline{\PP}$$
(see sections \ref{S-Prim}).

In particular, for any $n$, the sequence
$$Q(n)\mapright i \overline\PP
(n)\mapright{\overline\Delta}\overline{\PP}{\hat\otimes}{\overline \PP} (n)$$
is a left exact sequence of right $S_n$-modules. Recall besides that, for any finite group and 
any field $\kk$ of characteristic 0, every $\kk [G]$-module is projective. In particular, for any left
$S_n$-module, the tensor product $-\otimes_{S_n}M$
is an exact functor (see e.g. \cite[Sect. I.8]{Brown94})
and we have finally, for any $\S$-module $V$, a left exact sequence:
$$Q\circ V\rightarrow \overline\PP\circ V \mapright{\overline\Delta\circ
  V} (\overline{\PP} \hat\otimes {\overline\PP})\circ V$$
From (\ref{overline}) one has  ${\overline\PP}\circ V=\overline{\PP\circ V}$
and from (\ref{R-iso}) one has
$(\overline{\PP}{\hat\otimes} {\overline\PP})\circ
  V=\overline{\PP\circ
  V}\hat\otimes \overline{\PP\circ V}.$ As a consequence, $Q\circ V=\Prim(\PP\circ V)$.\hfill$\Box$

%%%%%%%%%%%%%%%Theoreme de cartier-Milnor-Moore

\section{A Cartier-Milnor-Moore theorem}

\subsection{Multiplicative Hopf operads}

Let $\PP$ be a connected Hopf operad. Then, $\PP$ has naturally the structure of a twisted 
Hopf $\PP$-algebra.
That is, there is a coproduct map from $\PP$ to $\PP\hat\otimes \PP$ which is a morphism of 
twisted $\PP$-algebras.

Assume that $\phi:U\rto\PP$ is a morphism of connected Hopf operad. In view of 
theorem \ref{T-hopfhatotimes}, this requirement amounts to the following condition. 
The morphism $\phi$, as any morphism of operads from $U$ to $\PP$, provides an arbitrary 
$\PP$-algebra with the structure of a $U$-algebra. In particular, the tensor product 
$A\hat\otimes B$ of two $\PP$-algebras, which is a $\PP$-algebra (since $\PP$ is a Hopf operad) 
carries naturally the structure of a $U$-algebra.

On the other hand, $\phi$ induces on $A$ and $B$ a structure of $U$-algebra and, since $U$ is 
a connected Hopf operad, the tensor product $A\hat\otimes B$ carries the structure of a $U$-algebra. 
The hypothesis that $\phi$ is a morphism of connected Hopf operads ensures that the two structures of 
$U$-algebras on $A\hat\otimes B$ are identical.

It follows in particular that the coproduct map from $\PP$ to $\PP\hat\otimes \PP$ is also, 
by restriction, a morphism of $U$-algebras, and
$\PP$ inherits from this construction the structure of a twisted Hopf $U$-algebra. More generally, 
we have:

\subsubsection{Proposition}\label{P-subHopf}\it 
Let $\phi: U\rto\PP$ be a morphism of connected Hopf operad.
Then, $\PP$ and, more generally, any twisted Hopf $\PP$-algebra, is naturally provided with
the structure of a twisted Hopf $U$-algebra.\rm

\subsubsection{Definiton}A {\sl multiplicative Hopf operad} $\PP$ is a connected Hopf 
operad together with a
morphism of connected Hopf operads $\phi:\As\rto\PP$.

Notice that the map $\As\rto\PP$ induces
a structure of twisted Hopf algebra on $\PP$ in view of theorem \ref{T-hopfhopf} 
and proposition \ref{P-subHopf},
since $\PP$ is a twisted connected Hopf $\PP$-algebra.

\subsection{Theorem}\label{T-CMM}\it Any multiplicative co-commutative Hopf operad is --as a twisted Hopf algebra-- 
the twisted enveloping algebra of its primitive elements.\rm

\medskip

\noindent{\sl Proof.} Let $\PP$ be a multiplicative co-commutative Hopf operad.
Recall first that the notion of enveloping algebra holds in the category of $\S$-modules -that is, 
a twisted (associative) algebra is naturally associated to any twisted Lie algebra (with the usual 
universal properties of enveloping algebras, see e.g. \cite{Joyal86}).
Recall also that the Cartier-Milnor-Moore
theorem holds for twisted connected Hopf algebras in any characteristic. More precisely,
if $H$ is a connected cocommutative twisted Hopf
algebra, $\Prim(H)$ carries naturally the structure of a twisted Lie algebra (this follows e.g.  from our theorem \ref{Prim}). The embedding $\Prim(H)\hookrightarrow H$ of the primitive
elements of $H$ into $H$ induces an isomorphism of twisted Hopf algebras
between the twisted enveloping algebra of $\Prim(H)$ and $H$. This was proven by Stover \cite{Stover93}; alternatively, as has been
pointed out by Fresse \cite[Appendix A]{Fr98}, the combinatorial proof of the
classical Cartier-Milnor-Moore given in \cite{Patras94} holds in any
graded linear symmetric monoidal category, and therefore applies to twisted Hopf algebras, which are
Hopf algebras in the category of $\S$-modules.
Our theorem follows.\hfill$\Box$

\section{The poisson operad}
In this section we assume that $\kk$ is of characteristic $0$.

Recall that a Poisson algebra $A$ is a commutative algebra with a unit $1$ provided with a Lie bracket $[,]$ 
which is a biderivation. That is, we have, besides the antisymmetry and Jacobi identities for $[,]$, 
the Poisson distributivity formula:
\begin{equation}\label{DistribPoisson}
[f,gh]=[f,g]h+g[f,h]
\end{equation}
In particular, we have $[f,1]=0$.

The simplest way to describe the Poisson operad $\Pois$ is through the corresponding functor:
\begin{equation}\label{PoissonOperad}
\Pois (V)=\Com\circ \Lie (V)
\end{equation}
where $\Com$ is the operad of commutative algebras with a unit, and $\Lie$ the Lie operad.
In concrete terms, an element of the free Poisson algebra over $V$ is a commutative polynomial 
in the Lie polynomials (the elements of $\Lie (V)$), and the bracket of two such commutative polynomials is
computed using (iteratively) the Poisson distributivity formula and the Lie bracket in $\Lie (V)$.

Due to the Poincar\'e-Birkhoff-Witt theorem, which states that $\As (V)$ is isomorphic to 
$\Com\circ \Lie (V)$ as analytic functors, $\Pois (V)$ and $\As (V)$ are isomorphic as analytic functors.
Therefore they are also isomorphic as $\S$-modules, as a consequence of the correspondance between
polynomial functors and symmetric group representations. We refer to  \cite[Appendix A]{McD95} for further details on analytic functors, polynomial functors and symmetric group representations. In particular,
$\Pois (n)$ is isomorphic, as a right $S_n$-module to the regular representation of $S_n$. 

Recall however that the Poincar\'e-Birkhoff-Witt theorem also holds in the category of $\S$-modules: the
enveloping algebra of a connected twisted Lie algebra $L$ is isomorphic, as a $\S$-module, to the 
free twisted commutative algebra over $L$ \cite[theorem 2]{Joyal86}.

We are going to show that the (classical) Poincar\'e-Birkhoff-Witt isomorphism between 
$\Pois (V)$ and $\As (V)$
can be lifted to the Hopf operadic setting, and understood directly by means of the Joyal's 
Poincar\'e-Birkhoff-Witt theorem for the twisted enveloping algebras of twisted Lie algebras.

\subsection{Primitive elements of the Poisson operad}
Recall that the Poisson operad is a connected cocommutative Hopf operad. Let $[,]$ and $\mu$ be 
the two generators
in $\Pois(2)$ representing the Lie structure and commutative structure. Then one has the following

$$\begin{cases} \mu|_{\emptyset}=1_0,& \\
\mu|_{S}=1_1,&\mbox{\rm for}\ |S|=1\end{cases}\quad  \mbox{\rm and}\quad
[,]|_T =0,\ \mbox{\rm for}\ |T|<2.$$
and
$$\delta(\mu)=\mu\otimes\mu, \quad
\delta([,])=[,]\otimes\mu+\mu\otimes [,].$$
The last two equations mean in terms of Poisson algebras, that
if $A$ and $B$ are two Poisson algebras, the tensor product $A\otimes
B$ is provided with a Poisson structure as follows. As a commutative
algebra, $A\otimes B$ is provided with the usual commutative product,
$(a_1\otimes b_1)\cdot (a_2\otimes b_2)=a_1a_2\otimes b_1b_2$.
The bracket on $A\otimes B$ is defined by: $[a_1\otimes b_1,a_2\otimes b_2]=a_1a_2\otimes [b_1,b_2]+[a_1,a_2]\otimes b_1b_2$.

\medskip

It follows in particular from the definition of $\delta$ that the natural inclusion gives
a morphism of connected Hopf operad
$\Com\rto\Pois$, so that the
proposition \ref{P-subHopf} applies: any twisted Hopf $\Pois$-algebra is
naturally provided with the structure of a commutative twisted Hopf algebra.

\subsection{Lemma}\it  The Lie operad is a suboperad of the operad of primitive elements in $\Pois$.\rm

\medskip

{\sl Proof.} Since the Lie operad is a sub-operad of $\Pois$ and is generated by the Lie bracket $[,]$, it is
enough, to prove the lemma, to check that $[,]\in \Pois (2)$ is a
primitive element. Thanks to theorem \ref{T-hopfhopf} one has to compute
$$\Delta([,])=\sum_{(1),(2)}\sum_{S\sqcup T=[2]}([,]_{(1)})|_S\otimes ([,]_{(2)})|_T.$$
Since $[,]|_{S}=0$ for $|S|<2$ the latter equality writes
$$\Delta([,])=1_0\otimes [,]+[,]\otimes 1_0,$$
and $[,]$ is primitive.\hfill $\Box$

\subsection{Theorem}\it The Poisson operad is naturally provided with
the structure of a commutative and cocommutative twisted Hopf
algebra. 
The suboperad of primitive elements of $\Pois$ is the Lie
operad. Moreover, as a twisted Hopf algebra, $\Pois$ is isomorphic to
the free 
commutative twisted algebra over the $\S$-module $\Lie$.\rm

\medskip

\noindent{\sl Proof.} As pointed out before, the inclusion $\Com\subset\Pois$ induces a 
morphism of connected Hopf operads
$\Com\rto\Pois$, thus, by composition with $\As \to\Com$, a morphism
of connected Hopf operads $\As\rto\Pois$. Then we can apply the theorem \ref{T-CMM} to $\Pois$.
Since $\Pois$ is also twisted commutative, the twisted Lie algebra structure on $\Prim(\Pois)$ is trivial 
(see \cite{Stover93}),
and $\Pois$ is isomorphic, as a Hopf algebra, to the free twisted commutative algebra over 
the $\S$-module of its primitive elements (recall that a free twisted $\PP$-algebra over a
 connected Hopf operad $\PP$ is naturally provided with a twisted Hopf
 $\PP$-algebra structure, so that a free commutative
 twisted algebra is naturally provided with a twisted commutative Hopf algebra structure).

It remains to prove that $\Lie$ is the set of primitive elements in
$\Pois$. We conclude by a dimensionnality argument based on the remark
that, 
if the $\S$-module $A$ is a sub $\S$-module of $B$, and if $dim_{\kk}
A(n)<\infty$ for all $n$, then, if the free twisted commutative 
algebras over $A$ and $B$ have the same dimension over $\kk$ in each degree, then $A=B$.

Recall from \cite[Prop.17]{PatReu04} that $\As$ is, as a twisted Hopf algebra, 
the enveloping algebra of $\Lie$. Due to Joyal's Poincar\'e-Birkhoff-Witt theorem, 
it follows that 
the dimension of the component of degree $n$ of the free commutative twisted algebra over 
$\Lie$ is equal to the 
dimension of $\As (n)$, that is, to $n!$. Besides, $\Pois (n)$ is isomorphic to the regular 
representation of $S_n$ 
as a $S_n$-module and, in particular, has dimension $n!$ as a vector space. Since $\Lie$ is 
contained in $\Prim(\Pois )$, 
and since, according to our previous arguments, the dimensions of the graded components of the 
free twisted commutative 
algebras over $\Lie$ and over $\Prim(\Pois )$ are equal, the theorem follows: 
$\Lie =\Prim(\Pois )$.\hfill$\Box$

%\begin{thebibliography}{}
%\end{thebibliography}

%\vspace{-0.3cm}

\bibliographystyle{amsalpha} 
\bibliography{bibliojuin06.bib}

\end{document}